\documentclass[12pt]{article}

\usepackage{amsfonts}

\begin{document}
\begin{center}
{\LARGE Integer Triangles with $R/r=N$}

\medskip

Allan MacLeod - University of Paisley

\medskip
(e-mail: allan.macleod@paisley.ac.uk)

\end{center}

\bigskip

{\small

{\bf Abstract}: We consider the problem of finding integer-sided
triangles with $R/r$ an integer, where $R$ and $r$ are the radii
of the circumcircle and incircle respectively. We show that such
triangles are relatively rare.

 }

\bigskip
\section{Introduction}
Let $FGH$ be a triangle with sides of integer length $f,g,h$.
Possibly the two most fundamental circles associated with the
triangle are the circumcircle which passes through the $3$
vertices, and the the incircle which has the $3$ sides as
tangents. The radii of these circles are denoted $R$ and $r$
respectively. Clearly $R>r>0$, and it is interesting to speculate
whether the (dimensionless) ratio $R/r$ could ever be an integer.

If we denote the distance between the centres of these two
circles as $d$, then it is a standard result in triangle geometry
that $d^2=R(R-2r)$, so that $R/r \ge 2$. Equality occurs for
equilateral triangles with $f=g=h=1$.

Basic trigonometry gives the formulae
\begin{equation}
R=\frac{f g h}{4 \Delta}  \hspace{2cm} r = \frac{\Delta}{s}
\end{equation}
with $s$ being the semi-perimeter $(f+g+h)/2$ and $\Delta$ the
area.

Thus
\begin{equation}
\frac{R}{r}=\frac{f g h s}{4 \Delta^2}
\end{equation}

Now, $\Delta=\sqrt{s(s-f)(s-g)(s-h)}$, so that
\begin{equation}
\frac{R}{r}=\frac{2 f g h}{(f+g-h)(f+h-g)(g+h-f)}
\end{equation}
so that, if $R/r=N$, with $N$ a strictly positive integer, we
must find integers $f,g,h$ such that
\begin{equation}
\frac{2 f g h}{(f+g-h)(f+h-g)(g+h-f)}=N
\end{equation}
which bears a very strong resemblance to the integer
representation problems in \cite{bgk} and \cite{bg}. In all cases
we look to express $N$ as a ratio of two homogeneous cubics in
$3$ variables.

Expressing equation (4) as a single fraction, we derive the cubic
\begin{equation}
N f^3   - N(g+h)f^2 - (g^2 N - 2 g h (N+1)+h^2 N)f +
N(g+h)(g-h)^2=0
\end{equation}

We cannot get very far with this form, but we can proceed quickly
if we replace $h$ by $2s-f-g$. This gives
\begin{equation}
2(4 N s - g(4N+1))f^2 - 2(g^2(4N+1)-2g s(6N+1)+8 N s^2)f + 8 N s
(g-s)^2=0
\end{equation}
and it should be noted that we only need to consider rational
$f,g,h$ and scale to integer values, since this would not affect
the value of $R/r$.

Since we want rational values for $f$, this quadratic must have a
discriminant which is a rational square, so that there must be a
rational $d$ with
\begin{equation}
d^2=(4N+1)^2g^4 - 4(2N+1)(4N+1)g^3s + 4(4N^2+8N+1)g^2s^2 - 16 N g
s^3
\end{equation}

Define $d=s^2y/(4N+1)$ and $g=s x/(4N+1)$, giving the quartic
\begin{equation}
y^2=x^4-4(2N+1)x^3+4(4N^2+8N+1)x^2-16 N x(4N+1)
\end{equation}
which can be transformed to an equivalent elliptic curve by a
birational transformation.

We find the curve to be
\begin{equation}
E_N: \; v^2=u^3+2(2N^2-2N-1)u^2+(4N+1)u
\end{equation}
with the transformation
\begin{equation}
\frac{g}{s}= \frac{v-(4N+1-(2N+1)u)}{(u-1)(4N+1)}
\end{equation}

As an example, consider the case of $N=7$, so that $E_7$ is the
curve $v^2=u^3+166u^2+29u$. It is moderately easy to find the
rational point $u=29/169,v=6902/2197$ which lies on the curve.
This gives $g/s=63/65$, and the equation for $f$ is
$-14f^2+938f+14560=0$, from which we find the representation
$f=-13,g=63,h=80$, clearly not giving a real-world triangle.

\section{Elliptic Curve Properties}
We have shown that integer solutions to equation (4) are related
to rational points on the curves $E_N$ defined in equation (9).
The problem is that equation (4) can be satisfied by integers
which could be negative as in the representation problems of
\cite{bgk} and \cite{bg}.

To find {\bf triangles} for the original form of the problem, we
must enforce an extra positivity constraint on $f,g,h$. To
investigate the effect of this, we must examine the properties of
the curves. We will try to keep this as simple as possible.

We first note that the discriminant of the curve is given by
\begin{displaymath}
\Delta = 256 N^3 (N-2) (4N+1)^2
\end{displaymath}
so that the curve is singular for $N=2$, which we now exclude
from possible values, having seen before that the equilateral
triangle gives $N=2$.

Given the standard method of addition of rational points on an
elliptic curve, see \cite{st}, the set of rational points forms a
finitely-generated group.

The points of finite-order are called torsion points, and we look
for these first. The point at infinity is considered the identity
of the group. The form of $E_N$ implies, by the Nagell-Lutz
theorem, that the coordinates of the torsion points are integers.

The points of order 2 are integer solutions of
\begin{displaymath}
u^3+2(2N^2-2N-1)u^2+(4N+1)u=0
\end{displaymath}
and it is easy to see that the roots are $u=0$ and $u=2N+1-2N^2
\pm \sqrt{N(N-2)}$. The latter two are clearly irrational and
negative for $N$ a positive integer, so there is only one point of
order 2.

The fact that there are $3$ real roots implies that the curve
consists of two components - an infinite component for $u \ge 0$
and a closed finite component (usually called the "egg").

Points of order $2$ allow one to look for points of order $4$,
since if $P=(j,k)$ has order $4$, $2P$ must have order $2$. For a
curve of the form given by (9), the u-coordinate of $2P$ must be
of the form $(j^2-4N-1)^2/4k^2$. Thus, we must have
$j^2=4N+1=(2r+1)^2$, so $N=r^2+r$. Substituting these values into
(9), we see that we get a rational point only if $r(r+2)$ is an
integer square, which never happens. There are thus no points of
order 4.

Points of order $3$ are points of inflexion of the curve. Simple
analysis shows that there are points of inflexion at $(1,\pm
2N)$. We can also use the doubling formula to show that there are
$2$ points of order $6$, namely $(4N+1,\pm 2N(4N+1))$.

The presence of points of orders $2,3,6$, together with Mazur's
theorem on possible torsion structures, shows that the torsion
subgroup must be isomorphic to $\mathbb{Z}6$ or $\mathbb{Z}12$.

The latter possibility would need a point $P$ of order $12$, with
$2P$ of order $6$, and thus an integer solution of
\begin{displaymath}
\frac{(j^2-(4N+1))^2}{4k^2}=4N+1
\end{displaymath}
implying that $N=t^2+t$. Substituting into this equation, we get
an integer solution if either $t^2-1$ or $t(t+2)$ are integer
squares - which they are not, unless $N=2$ which has been
excluded previously.

Thus the torsion subgroup is isomorphic to $\mathbb{Z}6$, with
finite points $(0,0)$, $(1,\pm2N)$, $(4N+1,\pm2N(4N+1))$.
Substituting these points into the $g/s$ transformation formula
leads to $g/s$ being $0$, $1$ or undefined, none of which lead to
a practical solution of the problem.

Thus, we must look at the second type of rational point - those
of infinite order.

\section{Practical Solutions}
If there are rational points of infinite order, Mordell's theorem
implies that there are $r$ generator points $G_1,\ldots,G_r$, such
that any rational point $P$ can be written
\begin{equation}
P= T+n_1G_1+\ldots+n_rG_r
\end{equation}
with $T$ one of the torsion points, and $n_1,\ldots,n_r$
integers. The value $r$ is called the {\bf rank} of the curve.

Unfortunately, there is no simple method of determining firstly
the rank, and then the generators. We used a computational
approach to estimate the rank using the Birch and Swinnerton-Dyer
conjecture. A useful summary of the computations involved can be
found in the paper of Silverman \cite{silv}.

Applying the calculations to a range of values of $N$, we find
several examples of curves with rank greater than zero, mostly
with rank $1$. A useful by-product of the calculations in the
rank $1$ case is an estimate of the height of the generator,
where the height gives an indication of number of digits in the
rational coordinates. For curves of rank greater than $1$ and
rank $1$ curves with small height, we can reasonably easily find
generators. However, when we backtrack the calculations to
solutions of the original problem, we hit a significant problem.

The elliptic curve generators all give solutions to equation (4),
but for the vast majority of $N$ values, these include at least
one negative value of $f,g,h$. Thus we find extreme difficulty in
finding real-life triangles with strictly positive sides. In
fact, for $3 \le N \le 99$, there are only $2$ values of $N$
where this occurs, at $N=26$ and $N=74$.

To investigate this problem, consider the quadratic equation (6),
but written as
\begin{equation}
f^2+(g-2s)f+\frac{4 N s (g-s)^2}{4 N s-g(4N+1)}
\end{equation}

The sum of the roots of this is clearly $2s-g$, but since
$f+g+h=2s$, this means that the roots of this quadratic are in
fact $f$ and $h$. Positive triangles thus require $s>0$, $g>0$,
$2s-g>0$ and $4Ns-g(4N+1)>0$, all of which reduce to
\begin{equation}
0 < \frac{g}{s} < \frac{4N}{4N+1}
\end{equation}

Looking at equation (10), we see that the analysis splits first
according as $u>1$ or $u<1$. Consider first $u>1$, so that, for
$g/s>0$ we need $v > 4N+1-(2N+1)u$. The line $v=4N+1-(2N+1)u$
meets $E_N$ in only two points, $(1,2N)$ and $(4N+1,-2N(4N+1))$
with the line actually being a tangent to the curve at the latter
point. Thus, if $u>1$ we need $v>0$ to give points on the curve
with $g/s>0$.

For the second half of the inequality with $u>1$, we need
$v<1+(2N-1)u$. The line $v=1+(2N-1)u$ meets $E_N$ only at
$(1,2N)$, so none of the points with $u>1,v>0$ are satisfactory.
Thus to satisfy (13) we must look in the range $u<1$. Firstly, in
the interval $[0,1)$, we have $g/s>0$, since the numerator and
denominator of (10) are negative. For the second half, however,
we need $v>1+(2N-1)u$, but the previous analysis shows this
cannot happen.

Thus, the only possible way of achieving real-world triangles is
to have points on the egg component. From the previous analysis it
is clear that any point on the egg leads to $g/s>0$. For the other
part of (10), we must consider where the egg lies in relation to
the line $v=1+(2N-1)u$. Since the line only meets $E_N$ at $u=1$,
the entire egg either lies above or below the line. The line meets
the u-axis at $u=-1/(2N-1)$, and the extreme left-hand end of the
egg is at $u=2N+1-2N^2 - \sqrt{N(N-2)}$ which is less than $-1$
for $N \ge 3$. Thus the entire egg lies above the line so
$v>1+(2N-1)u$ and so (10) is satisfied.

Thus, we get a real-life triangle if and only if $(u,v)$ is a
rational point on $E_N$ with $u<0$.

Consider now the effect of the addition $P+T=Q$ where $P$ lies on
the egg and $T$ is one of the torsion points. All of the five
finite torsion points lie on the infinite component. Since $E_N$
is symmetrical about the u-axis, $P, T, -Q$ all lie on a straight
line, and since the egg is a closed convex region, simple
geometry implies that $-Q$ and hence $Q$ must lie on the egg.
Similarly, if $P$ lies on the infinite component then $Q$ must
also lie on the infinite component.

Geometry also shows that $2P$ must lie on the infinite component
irrespective of where $P$ lies.

\begin{center}
TABLE 1\\Integer sided triangles with $R/r=N$\\
\begin{tabular}{cccc}
$\;$&$\;$&$\;$&$\;$ \\
N&$f$&$g$&$h$\\
2 & 1 & 1& 1 \\
 26 &     11 &     39  &    49\\
  74  &    259&     475 &    729\\
 218  &   115&     5239 &   5341\\
 250  &   97 &     10051&   10125\\
 314  &  177487799 &  55017780825 &    55036428301\\
 386  &   1449346321141 &  2477091825117 &  3921344505997\\
 394  &   12017&   2356695 &    2365193\\
 458  &   395  &   100989  &    101251\\
 586  &   3809 &   18411 &  22201\\
 602  &   833  &   14703 &  15523\\
 634  &   10553413  &  1234267713   &   1243789375\\
 674  &   535  &   170471 &     170859\\
 746  &   47867463&    6738962807 &     6782043733\\
 778  &   1224233861981  & 91266858701995   &   92430153628659\\
 866  &   3025  &  5629  &  8649
\end{tabular}
\end{center}

This shows that if the generators of $E_N$ all lie on the
infinite component then there is {\bf NO} rational point on the
egg, and hence no real-life triangle.

We have, for $1 \le N \le 999$, found $16$ examples of integer
triangles, which are given in Table 1. There are probably more to
be found, but these almost certainly come from rank $1$ curves
with generators having a large height and therefore difficult to
find. I am not sure that the effort to find more examples is
worthwhile.

A close look at the values of $N$, shows that they all satisfy $N
\equiv 2 \pmod{8}$. Is this always true? If so,{\bf WHY?}.

\section{Nearly-equilateral Triangles}
As mentioned in the introduction, if we have an equilateral
triangle with side $1$ then $N=2$. This suggests investigating
how close we can get to $N=2$ with non-equilateral integer
triangles. We thus investigate
\begin{displaymath}
\frac{R}{r}=2+\frac{1}{M}
\end{displaymath}
with $M$ a positive integer.

We can proceed in an almost identical manner to before, so the
precise details are left out, but we use the same names for the
lengths and semi-perimeter. The problem is equivalent to finding
rational points on the elliptic curve
\begin{equation}
F_M: v^2 = u^3 + (6M^2+12M+4)u^2 +(9M^4+4M^3)u
\end{equation}
with
\begin{displaymath}
\frac{g}{s}= \frac{v -(9M^3+4M^2-(5M+2)u)}{(u-M^2)(9M+4)}
\end{displaymath}

The curves $F_M$ have an obvious point of order $2$ at $(0,0)$,
and can be shown to have points of order $3$ at
$(M^2,\pm2M^2(2M+1))$ and order $6$ at
$(9M^2+4M,\pm2M(2M+1)(9M+4))$. In general these are the only
torsion points, none of which lead to a practical solution.

For $M=2k^2+2k$, however, with $k>0$, the elliptic curve has $3$
points of order $2$, which lead to the isosceles triangles with
$f=2k,g=h=2k+1$. This shows that we can get as close to $N=2$ as
we like with an isosceles triangle. If we reverse the process and
start with an isosceles triangle, we can show that $M$ must be of
the form $2k^2+2k$.

The curves $F_M$ have two components, the infinite one and the
egg, and, as before, we can show that real-life triangles can
only come from rational points on the egg. Numerical experiments
show that these are much more common than for $E_N$.

As an example, the values for $M=89$ are

 $f \; = \; 1018 \; 8073 \; 7479 \; 43$

 $g \; = \; 1093 \; 7217 \; 9616 \; 73$

 $h \; = \; 1106 \; 5215 \; 5663 \; 04$

which gives a triangle with angles $55.16^{\circ}$,
$61.78^{\circ}$ and $63.06^{\circ}$.

\vspace{2cm}
\begin{flushleft}
{\bf Postal Address:}

Dr. Allan J. MacLeod,

Department of Mathematics and Statistics,

University of Paisley,

High St.,

Paisley,

SCOTLAND

PA1 2BE
\end{flushleft}

\end{document}